\documentclass[a4paper,12pt, reqno]{amsart}
\usepackage{a4wide}
\usepackage{amssymb,amsfonts}
\usepackage[all,arc]{xy}
\usepackage{enumerate}
\usepackage{mathrsfs}
\usepackage[utf8]{inputenc}
\usepackage{amsmath,mathtools,tikz,graphicx,caption,subcaption}
\usepackage{xcolor}
\usepackage{amssymb}
\newtheorem{thm}{Theorem}[section]

\newtheorem{prop}[thm]{Proposition}
\newtheorem{lem}[thm]{Lemma}

\theoremstyle{definition}
\newtheorem{defn}[thm]{Definition}

\theoremstyle{remark}

\newtheorem*{acknowledgements*}{Acknowledgements}

\makeatletter
\def\blfootnote{\gdef\@thefnmark{}\@footnotetext}
\makeatother

\def\house#1{\setbox1=\hbox{$\,#1\,$}%
\dimen1=\ht1 \advance\dimen1 by 2pt \dimen2=\dp1 \advance\dimen2 by 2pt
\setbox1=\hbox{\vrule height\dimen1 depth\dimen2\box1\vrule}%
\setbox1=\vbox{\hrule\box1}%
\advance\dimen1 by .4pt \ht1=\dimen1
\advance\dimen2 by .4pt \dp1=\dimen2 \box1\relax}

\begin{document}
\title{Biquadratic fields having a non-principal euclidean ideal class}

\author{Jaitra Chattopadhyay and Subramani Muthukrishnan}
\address[Jaitra Chattopadhyay and Subramani Muthukrishnan]{Harish-Chandra Research Institute, HBNI\\
Chhatnag Road, Jhunsi\\
Allahabad - 211019\\
INDIA}
\email[Jaitra Chattopadhyay]{jaitrachattopadhyay@hri.res.in}
\email[Subramani Muthukrishnan]{msubramani@hri.res.in}
\begin{abstract}
H. W. Lenstra \cite{lenstra} introduced the notion of an Euclidean ideal class, which is a generalization of norm-Euclidean ideals in number fields. Later, families of number fields of small degree were obtained with an Euclidean ideal class (for instance, in \cite{hester1} and \cite{cathy}). In this paper, we construct certain new families of biquadratic number fields having a non-principal Euclidean ideal class and this extends the previously known families given by H. Graves \cite{hester1} and C. Hsu \cite{cathy}. 
\end{abstract}
\maketitle

\section{Introduction}

\noindent
Let $K$ be an algebraic number field and  $\mathcal{O}_K$ be its ring of integers. Following the standard notations, let $d_K$, $Cl_K$, $h_K$ and $\mathcal{O}_K^*$ be the discriminant, the ideal class group, the class number and the multiplicative group of units of $K,$ respectively. For an integral ideal $I$ of $\mathcal{O}_K$, $Nm(I)$ denotes its norm and $[I]$ stands for its ideal class. The class group $Cl_K$ and its size $h_K,$ which are invariants associated with $K,$ play a central role in number theory. As a consequence, we value results that shed light on $Cl_K$ and $h_K$.

\medskip

\noindent
An integral domain $R$ is said to be  an \textit{Euclidean domain} if there exists a function $\phi: R \to \mathbb{N}\cup \{0\}$ such that for every $\alpha, \beta \not= 0 \in R,$ there exists $\gamma \in R$ such that $\phi(\alpha - \beta \gamma) < \phi(\alpha).$  One of the classical problems in algebraic number theory is to classify all number fields which are Euclidean with respect to the absolute value of the standard norm map. By multiplicatively extending the norm map $Nm$ from $\mathcal{O}_K$ to $K \setminus \{0\}$, the above definition can be written as follows: for every $x \in K \setminus \{0\}$ there exists $\alpha \in \mathcal{O}_K$ such that $N(x-\alpha) <1.$ In 1979, H.W. Lenstra \cite{lenstra} made a crucial observation and generalized the aforementioned definition. He defined an integral ideal $C$ to be \textit{norm-Euclidean} if for every $x \in K \setminus \{0\}$ there exists $\alpha \in C$ such that $N(x-\alpha) <1.$ 

\medskip
\noindent
Further, he generalized the notion of norm-Euclidean ideals as follows. 

\begin{defn}\label{defn1} \cite{lenstra}
Let $R$ be a Dedekind domain and $\mathbb{I}$ be the set of all fractional ideals containing $R$. A fractional ideal $C$ of $R$ is said to be an \textit{Euclidean ideal} if there exists a function $\psi : \mathbb{I} \rightarrow W$, for some well-ordered set $W$, such that for all $I \in \mathbb{I}$ and all $x \in I^{-1}C \setminus C$, there exists some $y \in C$ such that 
\begin{equation}
\psi((x-y)I^{-1}C) < \psi (I).
\end{equation}
\end{defn}

\medskip
\noindent
Lenstra proved that if $C$ is an Euclidean ideal, then every ideal in the ideal class $[C]$ is also Euclidean. Thus we can unambiguously define Euclidean ideal classes.

\medskip

\noindent
It is a well known fact that if $\mathcal{O}_K$ is Euclidean, then its ideal class group $Cl_K$ is trival \cite{marcus}. Likewise, if $\mathcal{O}_K$ has an Euclidean ideal class then $Cl_K$ is cyclic (Theorem 1.6, \cite{lenstra}). More precisely, if $\mathcal{O}_K$ has an Euclidean ideal class $[C]$ then $Cl_K$ is a cyclic group and $[C]$ generates $Cl_K$. However, the converse is not always true, for $\mathbb{Q}(\sqrt{-d})$ for $d = 19,23,24,31,35,39,40,43,47$ has no Euclidean ideal class even though $Cl_K$ is cyclic (Proposition 2.1, \cite{lenstra}). 

\medskip

\noindent
Thus to classify all number fields in which the converse holds true is an intriguing problem and was addressed by Lenstra \cite{lenstra}. He proved under the assumption of the generalized Riemann Hypothesis (GRH), that all number fields with cyclic class group, except for the imaginary quadratic fields, have an Euclidean ideal class. In other words, if $rank(\mathcal{O}_K^*)\geq 1,$ then $\mathcal{O}_K$ has an Euclidean ideal class if and only if $Cl_K$ is cyclic, assuming the GRH for number fields. 

\medskip

\medskip

\noindent
H. Graves proved a useful growth result for certain algebraic number fields, which removed the appeal to GRH \cite{hester2}. We state the precise statement below, which will be necessary in the course our proof.

\begin{thm}\label{lem1} \cite{hester2}
Suppose that $K$ is a number field such that $|\mathcal{O}_{K}^*|=\infty$, and that $C$ is a non-zero ideal of $\mathcal{O}_K$. If $[C]$ generates the class group of $K$ and 
\begin{equation*}
\left| \left\{ \textrm{prime ideal } \wp \subset \mathcal{O}_{K} : N(\wp)\leq x, [\wp]=[C], \pi_{\wp} \mbox{ is onto } \right\} \right|   \gg  \frac{x}{(\log x)^2},
\end{equation*}
where $\pi_{\wp}$ is the canonical map from $\mathcal{O}_K^*$ to $(\mathcal{O}_K/\wp)^*$, then $[C]$ is an Euclidean ideal class.
\end{thm}

\medskip

\noindent
Using this growth result Ram Murty and Graves removed GRH for large class of number fields \cite{hester-murty}. The precise statement is as follows.

\begin{thm}
Let $K$ be a number field that is Galois over $\mathbb{Q}$ and has cyclic class group. If its Hilbert class field $H(K)$ has abelian Galois group over $\mathbb{Q}$ and if rank $(\mathcal{O}_K^*) \geq 4$, then 
$$
Cl_K=\langle[C]\rangle \mbox { if and only if } [C] \mbox{ is an Euclidean ideal class }.
$$
\end{thm}

\noindent
Using her growth result, Graves \cite{hester1} constructed an explicit number field, viz. $\mathbb{Q}(\sqrt{2},\sqrt{35})$, having a non-principal Euclidean ideal class. In this paper, we generalize this work and construct a new family of quartic fields of the form $\mathbb{Q}(\sqrt{2},\sqrt{p_1p_2})$ which has a non-principal Euclidean ideal class. The precise statement is given at the end of this section. It is worthwhile to note that the quartic field of the form $\mathbb{Q}(\sqrt{2},\sqrt{p_1p_2})$ has unit rank $3$. 

\medskip

\noindent
In \cite{cathy}, Hsu explicitly constructed a family of biquadratic and cyclic quartic fields with unit rank $3$ having an Euclidean ideal class, described below. 
\begin{thm} \cite{cathy}
Suppose $K$ is a quartic field of the form $\mathbb{Q}(\sqrt{q},\sqrt{kr}).$   Then $K$ has a non-principal Euclidean ideal class whenever $h_K = 2$. Here the integers $q,k,r$ are all primes $\geq 29$ and are all congruent to $1$ $\pmod{4}.$ 
\end{thm}

\medskip

\noindent
In this article we shall provide new families of number fields each having an Euclidean ideal class. This produces a new class of number fields other than the ones constructed by H. Graves \cite{hester1} and C. Hsu \cite{cathy}. Our main theorems are stated as follows.

\medskip

\noindent

\begin{thm}\label{th3}
Let $K_1=\mathbb{Q}(\sqrt{q},\sqrt{kr})$, where $q \equiv 3 \pmod 4$ and $k, r \equiv 1 \pmod 4$ are prime numbers. Suppose that $h_{K_1} = 2$. Then $K_1$ has an Euclidean ideal class.
\end{thm}

\medskip

\noindent

\begin{thm}\label{lastth}
Let $p,q$ be two prime numbers with $p, q \equiv 1 \pmod 4$. Then the biquadratic field $K_2=\mathbb{Q}(\sqrt{2},\sqrt{pq})$ with $h_{K_2}=2$ has an Euclidean ideal class.
\end{thm}

\medskip

\noindent
It is efficacious to note at this point that combining Hsu's result with the Theorem \ref{th3} above, we get a much larger class of biquadratic fields having an Euclidean ideal class. More precisely,

\medskip 

\noindent

\begin{thm}\label{larger}
Let $K=\mathbb{Q}(\sqrt{p}, \sqrt{qr})$, where $p$ is any rational prime and $q, r \equiv 1 \pmod 4$ are prime numbers. If $h_K=2$, then $K$ has an Euclidean ideal class.
\end{thm}

\medskip

\noindent
This article is organized as follows. In section $2,$ we recall some standard preliminary results. Section $3$ contains the computation of conductor and Hilbert class fields of certain biquadratic fields. The proof of the main theorems are presented in section $4$. Finally, we list the class numbers of some of the biquadratic fields of our interest in section $5$ and the table provides numerous examples of such fields with class number two.

\medskip

\section{Preliminaries}

\noindent
Let $K/\mathbb{Q}$ be a finite Galois extension with Galois group $G$. If $G$ is abelian, then by the famous Kronecker-Weber theorem, $K \subseteq \mathbb{Q}(\zeta_m)$ for some positive integer $m$. The smallest such integer $m$ is called the {\it conductor} of $K$. The following well known Proposition, which will be occasionally used in the next section, provides the conductor of real quadratic fields.

\medskip

\noindent
\begin{prop}\label{conductor-real}
The real quadratic field $\mathbb{Q}(\sqrt{m})$ has conductor $m$ if $m \equiv 1 \pmod{4}$ and $4m$ if $m \equiv 2,3 \pmod{4}.$
\end{prop}

\medskip

\noindent
For our purposes of this article, we recall a couple of standard results from class field theory. For that, we start with the following proposition.

\medskip

\noindent
\begin{prop} \label{Artin symbol-1}\cite{cox}
Let $L/K$ be a Galois extension of number fields with Galois group $G$. Let $\mathfrak{p}$ be a non-zero prime ideal in $\mathcal{O}_K$ which is unramified in $L$ and let $\mathfrak{P}$ be a prime ideal in $\mathcal{O}_L$ lying above $\mathfrak{p}$. Then there is a unique element $\sigma \in G$ such that for all $\alpha \in \mathcal{O}_L$, we have $$\sigma(\alpha) \equiv \alpha^{|\mathcal{O}_K/\mathfrak{p}|} \pmod {\mathfrak{P}}.$$ 
\end{prop}

\medskip

\noindent
\begin{defn}\label{Artin symbol}
The unique element $\sigma$ in Proposition \ref{Artin symbol-1} is called the {\it Artin symbol} and is often denoted by $\left(\frac{L/K}{\mathfrak{P}}\right)$.
\end{defn}

\medskip

\noindent
{\bf Remark:} Suppose $\mathfrak{P_1}, \mathfrak{P}_2, \ldots , \mathfrak{P}_g$ be all the primes in $\mathcal{O}_L$ lying above $\mathfrak{p}$. For $\tau \in G$, we have, $$\left(\frac{L/K}{\tau(\mathfrak{P}_i)}\right)=\tau \left(\frac{L/K}{\mathfrak{P}_i}\right) \tau^{-1} \mbox{ for all } i.$$ Since $G$ acts transitively on $\{\mathfrak{P_1}, \mathfrak{P}_2, \ldots , \mathfrak{P}_g\}$, the set $\left\{\left(\frac{L/K}{\mathfrak{P}_i}\right) : i = 1,2, \ldots , g \right\}$ is a conjugacy class in $G$. Thus by $\left(\frac{L/K}{\mathfrak{p}}\right)$, we unambiguously mean this conjugacy class.

\medskip

\noindent
Next, we define the {\it Dirichlet density} of a set of prime ideals as follows.

\medskip

\noindent
\begin{defn} \cite{cox}\label{density definition}
Let $K$ be a number field and let $\mathcal{S}$ be a set of prime ideals in $\mathcal{O}_K$. The Dirichlet density of the set $\mathcal{S}$ is defined to be $$\delta(\mathcal{S})=\lim_{s \to 1^+}\frac{\displaystyle\sum_{\mathfrak{p}\in \mathcal{S}}N(\mathfrak{p})^{-s}}{-\log (s-1)},$$ provided the limit exists. Here, $N(\mathfrak{p})$ stands for the quantity $|\mathcal{O}_K/\mathfrak{p}|$.
\end{defn}

\medskip

\noindent
Now, we are in a position to state the Chebotarev density theorem, that goes as follows.

\begin{thm}\textbf{(Chebotarev density theorem)}\label{lem7}
Let $L/K$ be a finite Galois extension of number fields with Galois group $G$ and let $\mathcal{C}$ be a conjugacy class in $G$. Then the Dirichlet density of the set
$$
\left\{\splitfrac{\mbox{ prime ideal }}{\mathfrak{p} \textrm{ in } \mathcal{O}_{K}} : \mathfrak{p} \mbox{ is unramified in L and } \left(\frac{\mathfrak{p}}{L/K}\right)=\mathcal{C} \right\}
$$
exists and equals $\frac{|\mathcal{C}|}{[L:K]},$ where $\left(\frac{\mathfrak{p}}{L/K}\right)$ is the Artin symbol of $\mathfrak{p}$.

\end{thm}

\medskip

\noindent
Finally, we state the following consequence of the Artin reciprocity law. 
\begin{lem}\label{hilbertclassfield}
Let $K$ be a number field and let $H(K)$ be its {\it Hilbert Class field} i.e., $H(K)$ is the unique maximal, unramified, abelian extension of $K$. Then the Galois group of $H(K)$ over $K$ is isomorphic to $Cl_K.$
\end{lem}

\noindent
Alternatively, the Hilbert class field $H(K)$  of $K$ is the unique maximal, abelian extension of $K$ such that the principal prime ideals of $K$ split completely in $H(K)$.

\medskip

\section{Computations for conductors and Hilbert class fields}
\noindent
The construction of the Hilbert class field of a number field is almost a century old classical problem. We are interested in the explicit construction of Hilbert class fields for certain family of biquadratic fields as these computations are essential to the proof of our main theorems. As mentioned, the conductors for biquadratic fields of the form $K=\mathbb{Q}(\sqrt{q},\sqrt{kr}),$ where $q,k$ and $r$ are primes, can be obtained by the following elementary arguments.

\begin{lem}\label{hester-suggestion-1}
Let $L/K$ be an extension of number fields. Assume that $L/\mathbb{Q}$ is a Galois extension and it is abelian. Then $f(K)$ divides $f(L)$, where $f(K), f(L)$ stands for the conductors of $K$ and $L$, respectively.
\end{lem}

\begin{proof}
Let $G$ be the Galois group of the extension $L/\mathbb{Q}$. As $G$ is abelian, the fundamental theorem of Galois theory implies that $K/\mathbb{Q}$ is also an abelian extension. Thus the conductor of $K$ is well-defined. 

\medskip

\noindent
Since $K \subseteq L \subseteq \mathbb{Q}(\zeta_{f(L)})$ and $f(K)$ is the smallest positive integer such that $K \subseteq \mathbb{Q}(\zeta_{f(K)})$, we have $$K \subseteq \mathbb{Q}(\zeta_{f(K)}) \subseteq \mathbb{Q}(\zeta_{f(L)}).$$ Since $\zeta_{f(L)}$ is a primitive root of unity and $\mathbb{Q}(\zeta_{f(K)}) \subseteq \mathbb{Q}(\zeta_{f(L)})$, we get that $f(K)$ divides $f(L).$

\end{proof}

\medskip

\noindent
The next lemma, which will often be helpful to us, provides the conductor of the compositum of two number fields in terms of those of the constituent fields.

\begin{lem}\label{hester-suggestion-2}
Let $L$ be an abelian number field and let $K_1$, $K_2$ be two subfields of $L$ such that $L$ is the compositum of $K_1$ and $K_2$. Then $f(L)=lcm(f(K_1), f(K_2))$, where $f(K_1), f(K_2)$ and $f(L)$ are the conductors of $K_1$, $K_2$ and $L$, respectively.
\end{lem}

\begin{proof}
Consider the following tower of number field extensions.

\begin{center}
\begin{tikzpicture}[node distance = 1cm, auto]
      \node (Q) {$\mathbb{Q}$};
      \node (F') [above of=Q, left of=Q] {$K_1$};
      \node (F) [above of=Q, right of=Q] {$K_2$};
      \node (L) [above of=F, right of=F] {$\mathbb{Q}(\zeta_{f(K_2)})$};
      \node (M2) [above of=F', left of=F'] {$\mathbb{Q}(\zeta_{f(K_1)})$};
      \node (K) [above of=Q, node distance = 2cm] {$L$};
      \node (K') [above of=K] {$\mathbb{Q}(\zeta_{f(L)})$};
      \draw[-] (Q) to node {} (F');
      \draw[-] (Q) to node {} (F);
      \draw[-] (F') to node {} (K);
      \draw[-] (F) to node {} (K);
      \draw[-] (M2) to node {} (K');
      \draw[-] (L) to node {} (K');
      \draw[-] (F') to node {} (M2);
      \draw[-] (F) to node {} (L);
      \draw[-] (K) to node {} (K');
      \end{tikzpicture}
\end{center}

\begin{center}
	(Fig. 3.1)
\end{center}

\noindent
Since $K_1$ and $K_2$ are subfields of an abelian extension $L/\mathbb{Q}$, by Lemma \ref{hester-suggestion-1}, $f(K_1)$ divides $f(L)$ and  $f(K_2)$ divides $f(L)$. Therefore, $lcm(f(K_1), f(K_2))$ divides $f(L)$.

\medskip

\noindent
On the other hand, using the following facts that $$L=K_1K_2 \subseteq \mathbb{Q}(\zeta_{f(K_1)})\mathbb{Q}(\zeta_{f(K_2)})=\mathbb{Q}(\zeta_{lcm(f(K_1), f(K_2)})$$ along with the minimality of the conductor implies that $\mathbb{Q}(\zeta_{f(L)}) \subseteq \mathbb{Q}(\zeta_{lcm(f(K_1), f(K_2)}).$ Thus $f(L)$ divides $lcm(f(K_1), f(K_2))$ and hence $f(L)=lcm(f(K_1), f(K_2))$.

\end{proof}

\medskip

\begin{lem}\label{lem3} 
The conductor $f(K)$ of $K=\mathbb{Q}(\sqrt{q},\sqrt{kr})$ is $4qkr$, where $q, k, r$ are primes such that $q \equiv 3 \pmod 4$ and $k, r \equiv 1 \pmod 4$.
\end{lem}

\begin{proof}
Let $K_1=\mathbb{Q}(\sqrt{q})$ and $K_2=\mathbb{Q}(\sqrt{kr})$ be two quadratic fields. By Proposition \ref{conductor-real}, we have $f(K_1)=4q$ and $f(K_2)=kr$. As $K=K_1K_2$, we get by Lemma \ref{hester-suggestion-2}, $$f(K)=lcm(f(K_1), f(K_2))=4qkr.$$
\end{proof}

\medskip

\noindent
We now use Lemma \ref{lem3} to find the Hilbert class field of the biquadratic fields considered above.
\begin{lem}\label{lem5}
Let $K=\mathbb{Q}(\sqrt{q}, \sqrt{kr})$, where $q, k \mbox{ and } r$ are prime numbers satisfying $q \equiv 3 \pmod 4$ and $k, r \equiv 1 \pmod 4$. If $K$ has class number two, then the Hilbert class field $H(K)$ is $\mathbb{Q}(\sqrt{q},\sqrt{k},\sqrt{r})$.  
\end{lem}

\begin{proof}
Let $K'=\mathbb{Q}(\sqrt{q},\sqrt{k},\sqrt{r})$ be a number field, where $q,k,r$ are defined as above. By applying Lemma \eqref{lem3} twice, we get $f(K')=lcm(f(K_1K_2), f(K_3))=lcm(4qk, r)=4qkr.$

\medskip

\noindent
Now, we are ready to compute the Hilbert class field of $K$. Artin reciprocity law gives (by Lemma \eqref{hilbertclassfield}),  $Gal(H(K)/K) \simeq Cl_K.$

\noindent
By our assumption on the has class number of $K$, the Hilbert class field $H(K)$ is a quadratic extension of $K.$ Therefore, it is sufficient to prove that $K'$ is an unramified extension of $K$.

\medskip

\noindent
From the tower of number fields, we have 

$$\mathbb{Q} \subsetneq K \subsetneq K' \subsetneq \mathbb{Q}(\zeta_{4qkr}).$$ 

\noindent
Thus, the prime ideals in $K$ lying above $2, q, k$ and $r$ may ramify $K'$. In fact, we prove that all such prime ideals are unramified for $K'$ over $K$.

\medskip

\noindent
We first prove this for the rational prime $2$.
\noindent
Let $L_1 = \mathbb{Q}(\sqrt{k},\sqrt{r})$ and let $\mathfrak{p}_1$, $\mathfrak{p}_2$, $\mathfrak{p}_3$ and $\mathfrak{p}_4$ be primes in $L_1$, $K$, $K'$ and $\mathbb{Q}(\zeta_{4qkr})$ respectively, all lying above $2$ and the containments are indicated in the diagram below.

\begin{figure}[htbp]
\centering
    \begin{tikzpicture}[node distance = 1cm, auto]
      \node (Q) {$\mathbb{Q}$};
      \node (L_1) [above of=Q, right of=Q] {$L_1$};
      \node (K) [above of=Q] {$K$};
      \node (K') [above of=K] {$K'$};
      \node (H) [above of=K'] {$\mathbb{Q}(\zeta_{4qkr})$};
      \draw[-] (Q) to node {} (L_1);
      \draw[-] (Q) to node {} (K);
      \draw[-] (K') to node {} (H);
      \draw[-] (L_1) to node {} (H);
      \draw[-] (L_1) to node {} (K');
      \draw[-] (K) to node {} (K');
\end{tikzpicture}
\hspace{3cm}
    \begin{tikzpicture}
    [node distance = 1cm, auto]
      \node (Q) {$2$};
      \node (L_1) [above of=Q, right of=Q] {$\mathfrak{p}_1$};
      \node (K) [above of=Q] {$\mathfrak{p}_2$};
      \node (K') [above of=K] {$\mathfrak{p}_3$};
      \node (H) [above of=K'] {$\mathfrak{p}_4$};
      \draw[-] (Q) to node {} (L_1);
      \draw[-] (Q) to node {} (K);
      \draw[-] (K') to node {} (H);
      \draw[-] (L_1) to node {} (H);
      \draw[-] (L_1) to node {} (K');
      \draw[-] (K) to node {} (K');
\end{tikzpicture}
\begin{center}
	(Fig. 3.2)
\end{center}
  
\end{figure}

\medskip

\noindent
Since $k, r \equiv 1 \pmod 4$, the prime $2$ is unramified in $\mathbb{Q}(\sqrt{k})$ and $\mathbb{Q}(\sqrt{r})$ and so in the compositum $L_1$ i.e., $e(\mathfrak{p}_1|2)=1$. And, since $q \equiv 3 \pmod 4$, $2$ is ramified in $\mathbb{Q}(\sqrt{q}).$ As $\mathbb{Q}(\sqrt{q}) \subseteq K \subseteq K',$ we have $e(\mathfrak{p}_2 | 2) > 1$ and  $e(\mathfrak{p}_3 | 2) > 1$. 

\medskip

\noindent
From the multiplicativity of the ramification indices, we get, \begin{equation*}
1 < e(\mathfrak{p}_3 | 2)=e(\mathfrak{p}_3 | \mathfrak{p}_1)e(\mathfrak{p}_1 | 2) \leq 2\cdot 1.
\end{equation*}

\noindent
Thus, $e(\mathfrak{p}_3 | 2)=2$ and $e(\mathfrak{p}_2 | 2)=2$ and hence $e(\mathfrak{p}_3 | \mathfrak{p}_2)=\frac{e(\mathfrak{p}_3| 2)}{e(\mathfrak{p}_2 | 2)}=1$. That is, $\mathfrak{p}_3$ is unramified over $\mathfrak{p}_2.$

\medskip

\noindent
By a close examination of the arguments given above, we immediately see that just by replacing $2$ with $q$, we get the unramifiedness for the prime $q$. 

\medskip

\noindent
For the other two primes, viz, $k$ and $r$, the proof goes through ad verbatim. That is, we simply need to replace $L_1$ by $L_2=\mathbb{Q}(\sqrt{q},\sqrt{r})$ for the prime $k$ and $L_1$ by $L_3=\mathbb{Q}(\sqrt{k},\sqrt{q})$ for the prime $r$. Thus all the prime ideals in $K$ lying above the rational primes $q, k$ and $r$ are unramified in $K'$. Hence, the Hilbert class field $H(K)=K'=\mathbb{Q}(\sqrt{q},\sqrt{k},\sqrt{r}).$
\end{proof}

\noindent
Next, we turn our attention to compute the conductor and the Hilbert class field of a class of biquadratic fields which is slightly different from the one we discussed above. The proof of the following lemmas are similar, but still we present it again to make it complete.

\begin{lem}\label{2ndthm}
Let $K_2=\mathbb{Q}(\sqrt{2},\sqrt{pq}),$ where $p, q$ are rational primes with $p, q \equiv 1 \pmod 4$. Then the conductor $f(K_{2})$ of $K_2$ is $8pq$.
\end{lem}

\begin{proof}
It is an application of Lemma \eqref{lem3}.
\end{proof}

\begin{lem}\label{2ndthm1}
Let $K_2=\mathbb{Q}(\sqrt{2},\sqrt{pq}),$ be the biquadratic field as in the above lemma. If the class number $h_{K_2}$ of $K_2$ is $2$, then the Hilbert class field $H(K_2)$ of $K_2$ is $\mathbb{Q}(\sqrt{2}, \sqrt{p}, \sqrt{q})$.
\end{lem}

\begin{proof}
By close examination of the extension of number fields, the proof is similar to the proof of the Lemma \eqref{lem5}. But, it is crucial to note that the prime $2$ is ramified in $K_2.$ Therefore, we only sketch a proof of unramification of primes in $K_2$ lying above $2$.
Let $K''=\mathbb{Q}(\sqrt{2}, \sqrt{p}, \sqrt{q}).$ Using our standard arguments, we can prove that $K''$ has conductor $8pq$. Since $K_2$ has class number $2$ and $K''/K_2$ is a quadratic extension, the Hilbert class field of $K_2$ is equal to $K''$ provided that $K''/K_2$ is an unramified extension. \\
\noindent
As $K \subseteq K'' \subseteq \mathbb{Q}(\zeta_8)\mathbb{Q}(\zeta_p)\mathbb{Q}(\zeta_q)=\mathbb{Q}(\zeta_{8pq}),$ we only need to take care of the ramification of the primes $2,p$ and $q$. Since the argument is similar for all the primes exactly like in the previous Lemma \ref{2ndthm}, we just outline it for the prime $2$.  

\medskip

\noindent
The prime $2$ is unramified in $\mathbb{Q}(\sqrt{p}, \sqrt{q})$ and  it is ramified in $\mathbb{Q}(\sqrt{2})$. In other words, $e(\mathfrak{p}_2|2)=1$ and $e(\mathfrak{p}_1|2)=2$.
Therefore,
$$1<e(\mathfrak{p}_4|2)=e(\mathfrak{p}_4|\mathfrak{p}_2)\cdot e(\mathfrak{p}_2|2)=e(\mathfrak{p}_4|\mathfrak{p}_2)\leq 2$$

\noindent
and hence $e(\mathfrak{p}_4|2)=1.$ Finally, $e(\mathfrak{p}_4|\mathfrak{p}_3)=\frac{e(\mathfrak{p}_4|2)}{e(\mathfrak{p}_3|2)}=1,$ that is, $2$ is unramified in $K''$ over $K_2$.

\section{Proof of the main theorems}

\subsection{Proof of Theorem \eqref{th3}}
\noindent
The main strategy of our proof is to find an arithmetic progression modulo the conductor, suitable enough to apply the growth result as in the Theorem \eqref{lem1} by Graves. That is to show that the number of prime ideals $\wp$ belonging to the non-principal ideal class with norm less than or equal to $x$ and a unit generating the group $(\mathcal{O}_K/\wp)^*$ is $\gg x/\log^2 x.$ 

\medskip

\noindent
We first recall the following result as pointed in \cite{hester1}, which is crucial to fulfill our goal.

\begin{thm}\label{lem2} \cite{hester1}
	Let $K$ be a totally real number field with conductor $f(K)$ and let $\{e_1, e_2, e_3\}$ be a multiplicatively independent set contained in $\mathcal{O}_{K} ^*$. If $l=lcm(16,f(K))$, and if $gcd(u,l)=gcd(\frac{u-1}{2},l)=1$ for some integer $u$, then 
	\begin{align*}
	\left|\left\{ \substack{\textrm{primes } \subseteq \mathcal{O}_K \\ \textrm{of degree one}} : N({\wp}) \equiv u \ (mod \ l), N(\wp)\leq x, \langle -1,e_{i}\rangle \twoheadrightarrow (\mathcal{O}_{K}/{\wp})^* \right\}\right| \gg \frac{x}{(\log x)^2},
	\end{align*}
	for at least one $i$.
\end{thm}

\medskip

\noindent
Since we confine ourself to totally real quartic fields, by the previous theorem, we have enough number of prime ideals in the field satisfying the growth result \cite{hester1}, but we need them to be in the non-principal ideal class. Here is the place where the construction of the Hilbert class field is crucial.

\medskip

\noindent
\textbf{Proof of Theorem \eqref{th3}:} Let $K_1=\mathbb{Q}(\sqrt{q},\sqrt{kr})$ be a biquadratic field as in the Lemma \eqref{lem3} which is totally real with signature $(4,0)$. The field $K_1$ has the unit rank $3$, and therefore, there exists a multiplicatively independent set consisting of $3$ elements.

\medskip

\noindent
We first claim the following: 

\begin{equation}\label{condition1}
\text{\parbox{.85\textwidth}{There exists an integer $u$ such that whenever $p \equiv u \pmod{4qkr},$ $f(\mathfrak{p}|p)=1$ and $f(\wp | p)=2$ hold true. }}
\end{equation}

\medskip

\noindent
Here, $f(\mathfrak{p}|p)$ and $f(\wp | p)$ are the residual degrees of the primes $\mathfrak{p}$ and $\wp$ in $K_1$ and $H(K_1)$ respectively, both lying above $p$.

\medskip

\noindent
The proof of our claim goes as follows. Let $\displaystyle\left(\frac{p}{K_1/\mathbb{Q}}\right)$ and $\displaystyle\left(\frac{p}{H(K_1)/\mathbb{Q}}\right)$ be the Artin symbols of the prime $p$ in the number field $K_1$ and its Hilbert class field $H(K_1)$, respectively.

\medskip

\noindent
Consider the following two sets of prime numbers in $\mathbb{Z}$, namely,
\begin{center}
	$
	X_{K_1}=\displaystyle\left\{p : p \mbox{ is prime and } \displaystyle\left(\frac{p}{K_1/\mathbb{Q}}\right)=1 \right\},
	$
\end{center}
and 
\begin{center}
	$
	X_{H(K_1)}=\displaystyle\left\{p : p \mbox{ is prime and } \displaystyle\left(\frac{p}{H(K_1)/\mathbb{Q}}\right)=1 \right\}.
	$
\end{center}
The containment of the sets $K_1 \subseteq H(K_1),$ implies that a prime $p$ splits completely in $K_1$ whenever it splits completely in $H(K_1)$. Therefore, $X_{H(K_1)} \subseteq X_{K_1}.$

\medskip

\noindent
Since both the groups $Gal(K_1/\mathbb{Q})$ and $Gal(H(K_1)/\mathbb{Q})$ are abelian, the identity element itself constitutes a conjugacy class. So, by Theorem \eqref{lem7}, the Dirichlet densities of the sets $X_{K_1}$ and $X_{H(K_1)}$ are $\frac{1}{4}$ and $\frac{1}{8}$, respectively. Since $X_{H(K_1)} \subseteq X_{K_1}$, the Dirichlet density of the difference $X_{K_1} \setminus X_{H(K_1)}$ is $\frac{1}{4} - \frac{1}{8} = \frac{1}{8}$. In particular, the set $X_{K_1} \setminus X_{H(K_1)}$ is infinite and thus any of its element can serve our purpose. Understanding the set $X_{K_1} \setminus X_{H(K_1)}$ more closely will be convenient for us. We do this now.

\medskip

\noindent
First, we recall a result about splitting of primes in a number field from \cite{marcus}.

\begin{thm}\label{marcus1}\cite{marcus}
Let $K$ be a number field and let $L$ and $M$ be two extensions of $K$. Let $\mathfrak{p}$ be a prime ideal in the ring of integers of $K$. If $\mathfrak{p}$ splits completely in both $L$ and $M$, then $\mathfrak{p}$ splits completely in the compositum field $LM$.
\end{thm}

\noindent
Using the above theorem, we characterize all the prime numbers in the set $X_{K_1} \setminus X_{H(K_1)}$.

\medskip

\noindent
We note that, from Theorem \ref{marcus1}, a prime $p$ splits completely in $K_1$ if and only if it splits completely in $\mathbb{Q}(\sqrt{q})$ as well as in $\mathbb{Q}(\sqrt{kr})$. By the theory of decomposition of primes in quadratic fields, we have: 

\medskip

\noindent
\begin{equation}\label{five}
\text{\parbox{.90\textwidth}{A prime $p$ splits completely in $\mathbb{Q}(\sqrt{q})$  and in $\mathbb{Q}(\sqrt{kr})$ iff $\displaystyle\left(\frac{q}{p}\right)=1 \mbox{ and } \displaystyle\left(\frac{kr}{p}\right)=1.$
}}
\end{equation}

\medskip

\noindent
Similarly, $p \not\in X_{H(K_1)}$ if and only if $p$ does not split completely in at least one of the fields $\mathbb{Q}(\sqrt{q})$, $\mathbb{Q}(\sqrt{k})$ or $\mathbb{Q}(\sqrt{r})$. That is 

\begin{equation}\label{nine}
\displaystyle\left(\frac{q}{p}\right)=-1 \mbox{ or } \displaystyle\left(\frac{k}{p}\right)=-1 \mbox{ or } \displaystyle\left(\frac{r}{p}\right)=-1.
\end{equation}

\medskip

\noindent
Thus, combining the equations \eqref{five} and \eqref{nine}, we get

\begin{equation}\label{condition2}
\text{\parbox{.85\textwidth}{$X_{K_1} \setminus X_{H(K_1)}=\displaystyle\left\{p : p \mbox{ is prime and } \displaystyle\left(\frac{q}{p}\right)=1 \mbox{ and } \displaystyle\left(\frac{k}{p}\right)=\displaystyle\left(\frac{r}{p}\right)=-1 \right\}$}}
\end{equation}

\medskip

\noindent
We proceed for the proof now.

\medskip

\noindent
Let $l=lcm(16,f(K_1))=lcm(16,4qkr)=16qkr$. In order to apply Theorem \eqref{lem2}, we are required to find an integer $u$ satisfying the condition \eqref{condition1} along with the following:

\begin{equation}\label{eqn8}
(i) ~ gcd(u, 16qkr)=1,
\end{equation}
\begin{equation}\label{eqn9}
(ii) ~ gcd\left(\frac{u-1}{2},16qkr\right)=1.
\end{equation}

\medskip

\noindent
We note that the condition $(ii)$ above is equivalent to  the following simultaneous congruences:
$$
\left\{\begin{array}{ll}
u \not\equiv 1 \pmod q;\\
u \not\equiv 1 \pmod r;\\
u \not\equiv 1 \pmod k;\\
u \not\equiv 1 \pmod 4.
\end{array}\right. 
$$

\medskip

\noindent
In other words, it suffices to find a prime number $w \in X_{K_1} \setminus X_{H(K_1)}$ satisfying \eqref{condition2}, \eqref{eqn8} and \eqref{eqn9}. 

\medskip

\noindent
Here, we quote the following result due to P. Pollack \cite{pollack} which will be of our help to prove the existence of $w$.

\begin{thm}\cite{pollack}\label{residue} 
	Let $p \geq 5$ be a prime number. Then there exists a prime number $q<p$ such that $\displaystyle\left(\frac{q}{p}\right)=-1$ and $q \equiv 3 \pmod{4}$.
\end{thm}

\noindent
Now we are ready to find our desired $w$. By the aforementioned Theorem \eqref{residue}, there exist prime numbers $p_1, p_2$ and $p_3$ with $p_1 < q$, $p_2 < k$, $p_3 < r$ and satisfying 

\begin{equation}\label{ten}
\displaystyle \left(\frac{p_1}{q}\right)=\displaystyle \left(\frac{p_2}{k}\right)=\displaystyle \left(\frac{p_3}{r}\right)=-1 \mbox{ and } p_1, p_2, p_3 \equiv 3 \pmod 4.
\end{equation}

\noindent
Now, we look for a solution to the following system of simultaneous congruences.

\begin{equation}\label{eleven}
\left\{\begin{array}{ll}
x \equiv p_1 \pmod q;\\
x \equiv p_2 \pmod k;\\
x \equiv p_3 \pmod r;\\
x \equiv 3 \pmod 4.
\end{array}\right.
\end{equation}

\medskip

\noindent
The Chinese remainder theorem guarantees a unique solution modulo $4qkr$ to the above system congruences, say $x_0 $.

\medskip

\noindent
Since $gcd(x_0, 4qkr)=1$, by Dirichlet's Theorem for primes in arithmetic progression, there exist infinitely many prime numbers $l'$ satisfying $l' \equiv x_0 \pmod{4qkr}$. We pick such a prime number $l'$ and call it $w$. We note that $w$ satisfies conditions $(i)$ and $(ii)$ above. It remains to check that $w$ satisfies \eqref{condition2}.

\medskip

\noindent
For that, using the law of quadratic reciprocity along with the facts $q \equiv 3 \pmod 4, k,r \equiv 1 \pmod 4$ and $w \equiv 3 \pmod 4$, we get $$\displaystyle\left(\frac{q}{w}\right)=\displaystyle\left(\frac{w}{q}\right)(-1)^{\frac{w-1}{2}}=-\displaystyle\left(\frac{w}{q}\right)=-\displaystyle\left(\frac{p_1}{q}\right)=1.$$

\medskip

\noindent
Also, $$\displaystyle\left(\frac{k}{w}\right)=\displaystyle\left(\frac{w}{k}\right)=\displaystyle\left(\frac{p_2}{k}\right)=-1$$ and $$\displaystyle\left(\frac{r}{w}\right)=\displaystyle\left(\frac{w}{r}\right)=\displaystyle\left(\frac{p_3}{r}\right)=-1.$$

\medskip

\noindent
Thus, the prime number $w$ satisfies the conditions \eqref{condition2}, \eqref{eqn8} and \eqref{eqn9}. Hence, for any prime $p$ satisfying $p \equiv w \pmod{4qkr},$ we have that $f(\wp | \mathfrak{p})=2$. That is, $\mathfrak{p}$ does not split completeley in $H(K_1)$. Therefore, using lemma \ref{lem5}, we can say that $\mathfrak{p}$ is not a principal ideal and hence it will be a generator of $Cl_{K_1}$. This will prove that $K_1$ has an Euclidean ideal class.

\medskip

\noindent

\medskip

\noindent
\textbf{Proof of Theorem \eqref{lastth}}

\medskip

\noindent
The proof is almost similar to that of the previous one. By following the argument given above, \eqref{condition2} reads the following.

\begin{equation}\label{twelve}
\text{\parbox{.85\textwidth}{$X_{K_2} \setminus X_{H(K_2)}=\displaystyle\left\{s : s\mbox{ is prime and } \displaystyle\left(\frac{2}{s}\right)=1 \mbox{ and } \displaystyle\left(\frac{p}{s}\right)=\displaystyle\left(\frac{q}{s}\right)=-1 \right\}$}}
\end{equation}

\noindent
In this case we have, $l=lcm(16,f(K_2))=lcm(16,8pq)=16pq$. So, we need to find an integer $u$ such that the following two conditions hold: $(i)$ $gcd(u, 16pq)=1$ and $(ii)$ $gcd(\frac{u-1}{2},16pq)=1$. Again, condition $(ii)$ can be replaced by the following equivalent conditions : $ u \not\equiv 1 \pmod p;$ $u \not\equiv 1 \pmod q;$ and $u \not\equiv 1 \pmod 4$.

\medskip

\noindent
Again using the Theorem \eqref{residue}, we can choose prime numbers $p_1 \equiv 3 \pmod{4}$ and $p_2 \equiv 3 \pmod{4}$ such that $$p_1 < p, p_2 < q \mbox{ and } \displaystyle\left(\frac{p_1}{p}\right)=\displaystyle\left(\frac{p_2}{q}\right)=-1.$$

\noindent
Let $x_0$ be a unique solution $\pmod{8pq}$ to the system of congruences below. 
\begin{equation}\label{thirteen}
\left\{\begin{array}{ll}
x \equiv p_1 \pmod p;\\
x \equiv p_2 \pmod q;\\
x \equiv 7 \pmod 8.
\end{array}\right.
\end{equation}

\noindent
By Dirichlet's theorem for primes in an arithmetic progression, there exist infinitely many primes $w$ satisfying $w \equiv x_0 \pmod{8pq}$. We choose such a prime number $w$. Then $$\displaystyle\left(\frac{2}{w}\right)=1 \mbox{ since } w \equiv 7 \pmod 8.$$ 

\noindent
Also, by the quadratic reciprocity law, we get $$\displaystyle\left(\frac{p}{w}\right)=\displaystyle\left(\frac{w}{p}\right)=\displaystyle\left(\frac{x_0}{p}\right)=\displaystyle\left(\frac{p_1}{p}\right)=-1.$$ and $$\displaystyle\left(\frac{q}{w}\right)=\displaystyle\left(\frac{w}{q}\right)=\displaystyle\left(\frac{x_0}{q}\right)=\displaystyle\left(\frac{p_2}{p}\right)=-1.$$

\medskip

\noindent
We take $u=w$ and clearly it fulfills all our requirements. This completes the proof.

\end{proof}

\section{List of quartic fields and class numbers}

\noindent
There are numerous number of examples of quartic fields of the form $K=\mathbb{Q}(\sqrt{q},\sqrt{kr})$ with class number $2$. We list some quartic fields of our interest and its class number. The authors have computed the class number of such fields using Sage.
\begin{table}[h!]
	\centering
	\begin{tabular}{||c c c c c c c c c c||} 
		\hline
		$(q,k,r)$ & $h_K$ & $(q,k,r)$ & $h_K$ & $(q,k,r)$ & $h_K$ & $(q,k,r)$ & $h_K$ & $(q,k,r)$ & $h_K$ \\ [0.5ex] 
		\hline\hline
(3, 5, 13) & 2 & (3, 5, 17) & 2 & (3, 5, 29) & 8 & (3, 5, 37) & 2 & (3, 5, 41) & 4 \\
(3, 5, 53) & 4 & (3, 5, 61) & 4 & (3, 5, 73) & 4 & (3, 5, 89) & 4 & (3, 5, 97) & 4 \\
(3, 5, 101) & 8 & (3, 5, 109) & 4 & (3, 5, 113) & 2 & (3, 5, 137) & 6 & (3, 5, 149) & 4 \\
(3, 5, 157) & 6 & (3, 5, 173) & 4 & (3, 5, 181) &  8 & (3, 5, 193) & 4 & (3, 5, 197) & 12 \\
(3, 5, 229) & 4 & (3, 13, 5) & 2 & (3, 13, 17) & 4 & (3, 13, 29) & 4 & (3, 13, 37) & 4\\
(3, 13, 41) & 6 & (3, 13, 53) & 8 & (3, 13, 61) & 16 & (3, 13, 73) & 4 & (3, 13, 89) & 2\\
(3, 13, 97) & 4 & (3, 13, 101) & 4 & (3, 13, 109) & 8 & (3, 13, 113) & 4 & (3, 13, 137) & 2 \\
(3, 13, 149) & 6 & (3, 13, 157) & 8 & (3, 13, 173) & 4 & (3, 13, 181) & 8 & (3, 13, 193) & 4 \\ 
(3, 13, 197) & 2 & (3, 13, 229) & 8 & (3, 17, 5) & 2 & (3, 17, 13) & 4 & (3, 17, 29) & 2 \\
(3, 17, 37) & 2 & (3, 17, 41) & 12 & (3, 17, 53) & 12 &(3, 17, 61) & 2 & (3, 17, 73) & 4 \\
(3, 17, 89) & 4 & (3, 17, 97) & 4 & (3, 17, 101) & 4 & (3, 17, 109) & 2 & (3, 17, 113) & 4 \\ 
(3, 17, 137) & 4 & (3, 17, 149) & 8 & (3, 17, 157) & 8 & (3, 17, 173) & 6 & (3, 17, 181) & 2 \\
(3, 17, 193) & 12 & (3, 17, 197) & 2 & (3, 17, 229) & 4 & (3, 29, 5) & 8 & (3, 29, 13) & 4 \\
(3, 29, 17) & 2 & (3, 29, 37) & 2 & (3, 29, 41) & 2 & (3, 29, 53) & 4 & (3, 29, 61) & 2 \\
(3, 29, 73) & 4 & (3, 29, 89) & 6 & (3, 29, 97) & 4 & (3, 29, 101) & 4 & (3, 29, 109) & 4 \\
(3, 29, 113) & 2 & (3, 29, 137) & 6 & (3, 29, 149) & 20 & (3, 29, 157) & 2 & (3, 29, 173) & 6 \\ [1ex] 
		\hline
	\end{tabular}
\end{table}

\medskip

\noindent
Now, we list the class number of biquadratic fields of the form $\mathbb{Q}(\sqrt{2},\sqrt{pq}).$

\begin{table}[h!]
	\centering
	\begin{tabular}{||c c c c c c c c c c||} 
		\hline
		$(2,k,r)$ & $h_K$ & $(2,k,r)$ & $h_K$ & $(2,k,r)$ & $h_K$ & $(2,k,r)$ & $h_K$ & $(2,k,r)$ & $h_K$ \\ [0.5ex]
		\hline\hline
(2, 5, 13) & 4 & (2, 5, 17) & 2 & (2, 5, 29) & 4 & (2, 5, 37) & 2 & (2, 5, 41) & 4 \\
(2, 5, 53) & 4 & (2, 5, 61) & 2 & (2, 5, 73) & 6 & (2, 5, 89) & 4 & (2, 5, 97) & 2  \\
(2, 5, 101) & 4 & (2, 5, 109) & 6 & (2, 5, 113) & 4 & (2, 5, 137) & 4 & (2, 5, 149) & 2 \\
(2, 5, 157) & 6 & (2, 5, 173) & 2 & (2, 5, 181) & 4 & (2, 5, 193) & 2 & (2, 5, 197) & 12 \\
(2, 5, 229) & 4 & (2, 13, 5) & 4 & (2, 13, 17) & 4 & (2, 13, 29) & 2 & (2, 13, 37) & 2 \\ 
(2, 13, 41) & 4 & (2, 13, 53) & 4 & (2, 13, 61) & 4 & (2, 13, 73) & 2 & (2, 13, 89) & 2 \\
(2, 13, 97) & 2 & (2, 13, 101) & 8 & (2, 13, 109) & 2 & (2, 13, 113) & 4 & (2, 13, 137) & 4 \\ 
(2, 13, 149) & 6 & (2, 13, 157) & 2 & (2, 13, 173) & 4 & (2, 13, 181) & 6 & (2, 13, 193) & 2 \\
(2, 13, 197) & 2 & (2, 13, 229) & 12 & (2, 17, 5) & 2 & (2, 17, 13) & 4 & (2, 17, 29) & 2 \\
(2, 17, 37) & 2 & (2, 17, 41) & 12 & (2, 17, 53) & 4 & (2, 17, 61) & 2 & (2, 17, 73) & 4 \\
(2, 17, 89) & 8 & (2, 17, 97) & 8 & (2, 17, 101) & 4 & (2, 17, 109) & 6 & (2, 17, 113) & 4 \\
(2, 17, 137) & 8 & (2, 17, 149) & 8 & (2, 17, 157) & 8 & (2, 17, 173) & 6 & (2, 17, 181) & 2 \\
(2, 17, 193) & 24 & (2, 17, 197) & 2 & (2, 17, 229) & 16 & (2, 29, 5) & 4 & (2, 29, 13) & 2 \\
\hline
\end{tabular}
\end{table}

\begin{table}[h!]
	\centering
\begin{tabular}{||c c c c c c c c c c||} 
		\hline
		$(q,k,r)$ & $h_K$ & $(q,k,r)$ & $h_K$ & $(q,k,r)$ & $h_K$ & $(q,k,r)$ & $h_K$ & $(q,k,r)$ & $h_K$ \\ [0.5ex] 
		\hline\hline
(2, 29, 17) & 2 & (2, 29, 37) & 4 & (2, 29, 41) & 4 & (2, 29, 53) & 2 & (2, 29, 61) & 2 \\ 
(2, 29, 73) & 2 & (2, 29, 89) & 2 & (2, 29, 97) & 14 & (2, 29, 101) & 4 & (2, 29, 109) & 4\\
(2, 29, 113) & 4 & (2, 29, 137) & 12 & (2, 29, 149) & 10 & (2, 29, 157) & 2 & (2, 29, 173) & 2 \\
(2, 29, 181) & 8 & (2, 29, 193) & 2 & (2, 29, 197) & 16 & (2, 29, 229) & 4 & (2, 37, 5) & 2 \\
(2, 37, 13) & 2 & (2, 37, 17) & 2 & (2, 37, 29) & 4 & (2, 37, 41) & 4 & (2, 37, 53) & 2\\
(2, 37, 61) & 4 & (2, 37, 73) & 8 & (2, 37, 89) & 2 & (2, 37, 97) & 2 & (2, 37, 101) & 2 \\
		\hline
\end{tabular}
\end{table}

\pagebreak
\noindent{\bf Acknowledgements.} We express our heartfelt gratitude towards Prof. R. Thangadurai for his encouragement throughout this project and valuable comments to improve the quality of the paper. We sincerely thank Dr. Hester Graves for going through the manuscript meticulously several times and her suggestions which greatly improved the presentation of the paper. We are also thankful to Prof. K. Srinivas and Prof. T.R. Ramadas for their careful reading of an earlier version of this manuscript and for their useful comments. The second author would like to thank the Roman Number Theory Association and Prof. Francesco Pappalardi for their financial support to visit the Universit\`{a} Roma Tre where this work was initiated. We also thank Mr. R. Dixit for his valuable suggestions regarding the computational part. We are grateful to the Dept. of Atomic Energy, Govt. of India and Harish-Chandra Research Institute for providing financial support to carry out this research.

\end{document}